\documentclass{article}

\usepackage{amsmath, amssymb, enumerate, cite}

\usepackage[final]{graphicx}

\newtheorem{theorem}{Theorem}
\newtheorem{corollary}[theorem]{Corollary}
\newtheorem{remark}[theorem]{Remark}
\newtheorem{lemma}[theorem]{Lemma}
\newtheorem{definition}[theorem]{Definition}

\newcommand{\proof}{ {\sc Proof.\quad}}
\newcommand{\pend}{ \hfill $\square$ \\}

\numberwithin{equation}{section}  
\numberwithin{figure}{section}    
\numberwithin{table}{section}     
\numberwithin{theorem}{section}

\newcommand{\of}[1]{\ensuremath{\left( #1 \right)}}
\newcommand{\cb}[1]{\ensuremath{ \left\{ #1 \right\} }}
\newcommand{\sqb}[1]{\ensuremath{ \left[ #1 \right] }}
\newcommand{\st}{\,|\;}

\newcommand{\eps}{\ensuremath{\varepsilon}}

\renewcommand{\P}{\ensuremath{\mathcal{P}}}

\newcommand{\Min}{{\rm Min\,}}

\newcommand{\wmin}{{\rm wmin\,}}

\newcommand{\dom}{{\rm dom \,}}
\newcommand{\cl}{{\rm cl \,}}
\newcommand{\co}{{\rm co \,}}

\newcommand{\Int}{{\rm int\,}}

\usepackage[
colorlinks=true, linkcolor=blue, citecolor=blue, urlcolor=blue, filecolor=blue
]{hyperref}

\begin{document}
\title{Remarks on\\   Lagrange Multiplier Rules in Set Valued Optimization}

\author{Carola Schrage\footnote{Carola.Schrage@unibz.it}
}

\date{Received: date / Accepted: date}

\maketitle


\begin{abstract}
In this note, three Lagrange multiplier rules introduced in the literature for set valued optimization problems are compared. A generalization of all three results is given which proves that under rather mild assumptions, $x$ is a weak solution to the constrained problem, if and only if it is a weak solution to the Lagrangian.

{\bf Keywords:} Lagrange multiplier rule, Set valued optimization, weak efficiency 
\end{abstract}

\section{Introduction}
In \cite{DhingraLalitha2016}, an approximate weak lattice solution concept to set optimization problems is studied and a problem with constraints is studied. The authors formulate a Lagrangian function to this, using linear operators as dual variables and in Theorem 3.2 prove a Lagrange multiplier rule under rather strong assumptions on the image of the approximate solution $x_0$ of the primal problem.
In \cite[Theorem 4.1]{HernandezMarin2007lagrangian}, a different set of assumptions is applied to prove a Lagrange multiplier rule for weak lattice solutions. A different weak solution concept is studied in \cite{corley87}, also providing a Lagrange multiplier rule in Theorem 4.1.

In the present note, we will relax the assumptions  applied in these publications, and on the other also aim to strengthen the conclusions possible even under the more general assumptions. The resulting Lagrange multiplier rule will be a generalizations of \cite{corley87} and \cite{HernandezMarin2007lagrangian}, as well as of \cite{DhingraLalitha2016}. 

 We will point out that the convexity assumption on the primal function $f$ and the constraints is weakest in \cite[Theorem 3.2]{DhingraLalitha2016} and strongest in \cite[Theorem 4.1]{corley87}, while the remaining assumptions are weakest in \cite[Theorem 4.1]{corley87} and strongest in \cite[Theorem 3.2]{DhingraLalitha2016}. Surprisingly, the result in \cite[Theorem 4.1]{HernandezMarin2007lagrangian} is intermediate between both other results, especially, the assumptions in \cite[Theorem 3.2]{DhingraLalitha2016} are, but for the convexity assumption, stronger than those of \cite[Theorem 4.1]{HernandezMarin2007lagrangian}.
It is also notable that the assumptions in \cite{HernandezMarin2007lagrangian} and \cite{DhingraLalitha2016} are actually strong enough to guarantee that the weak lattice solution to the primal problem is in fact a weak solution with respect to the vector criteria.
 
In the next section, we will gather several definitions and provide results needed in the main sections. The third section collects the Lagrange multiplier rules from \cite{corley87}, \cite{HernandezMarin2007lagrangian} and \cite{DhingraLalitha2016} while in the fourth section, a generalization of all three results is presented and the connection between the known Multiplier Rules is established. A final section summarizes the results.

\section{Setting and Basic Results}

Throughout this note,
$X$ is a nonempty set, $Y$, $Z$ are Hausdorff topological vector spaces with duals $Y^*$, $Z^*$ and power sets $\P(Y)$ and $\P(Z)$ respectively, including $\emptyset$ and $Y$, or $Z$. The set $K\subseteq Y$ is a convex cone  with $0\in K$ and nonempty topological interior $\Int K$ and  likewise $C\subseteq Z$ is a convex cone with $0\in K$ and $\Int C\neq\emptyset$.
The set $\P(Y)$ is ordered by means of
\[
A\leq B\quad\Leftrightarrow B\subseteq A+K
\]
and $\P(Z)$ is ordered by means of
\[
A\leq B\quad\Leftrightarrow B\subseteq A+C,
\]
compare \cite{Kuroiwa1998natural} or \cite{FivePersonSoMF} on the order relations on power sets of linear spaces.

The (positive) dual cones of $K$ and $C$ are denoted as $K^+$ and $C^+$, respectively. If $e\in \Int K$, then $y^*(e)<0$ is true for all $y^*\in K^+\setminus\cb{0}$.
Throughout the text, $e\in\Int K$ is assumed and $\eps\geq 0$.

A set $A\subseteq Y$ is $K$-bounded, iff for every $0$-neighbourhood $U\subseteq Y$ there exists a $t>0$ such that $A\subseteq tU+K$. Especially, $A\subseteq Y$ is $K$--bounded, if an only if $A\subseteq -te+K$ is true for some $t>0$.

To any set valued function $f$, the domain of $f$ is $\dom f=\cb{x\in X\st f(x)\neq\emptyset}$.

A function $f:X\to\P(Y)$ is said to be convex, iff
\[
\forall x_1, x_2\in X,\, t\in\sqb{0,1}:\quad tf(x_1)+(1-t)f(x_2) \subseteq f(tx_1+(1-t)x_2)+K.
\]
especially, if $f$ is convex, then $f(x)+C$ is a convex set for all $x\in X$ and $\bigcup\limits_{x\in X}f(x)+C$ is convex, too.

\begin{lemma}\label{lem:convex}
Define
\[
Q=\cb{(y,z)\in Y\times Z\st \exists x\in X:\, y\in f(x)+K,\, z\in g(x)+C}.
\]
If $f$ and $g$ are convex functions, then $Q$ is a convex set. 
\end{lemma}
\proof
Let $f$, $g$ be convex functions, $(y_1,z_1)\in f(x_1)\times g(x_1)+K\times C$, $(y_2,z_2)\in f(x_2)\times g(x_2)+K\times C$ and $t\in\sqb{0,1}$, then
\[
t(y_1,z_1)+(1-t)(y_2,z_2)\in f(tx_1+(1-t)x_2)\times g(tx_1+(1-t)x_2)+K\times C\subseteq Q,
\]
proving the statement.
\pend

A set valued function $f:X\to\P(Y)$ is called (closely) convex like, if $\bigcup\limits_{x\in x}f(x)+K$ ($\cl\bigcup\limits_{x\in x}f(x)+K$) is convex, compare \cite[Definition 2.5]{frenk1999classes} and \cite[Definition 2.9]{Jahn04}. Thus  if $f$ and $g$ are convex, then $(f\times g):X\to\P(Y\times Z)$ is convex like and especially closely convex like.

\begin{lemma}\label{lemma3}\cite[Lemma 2.1]{CrespiRoccaSchrage14}
To any set $A\subseteq Y$, it holds $A+\Int K=\Int(A+K)$ and $A+K+\Int K=\Int (A+K)$.
\end{lemma}

In the following, we will investigate the connection between solutions of
\begin{align}\tag{$P$}\label{eq:P}
\Min f(x) \text{, such that $x\in M$}.
\end{align}
for $M=\cb{x\in X\st 0\in g(x)+Z}$ and those of
\begin{align}\tag{$LP_T$}\label{eq:LP_T}
\Min L(x,T) \text{, such that $x\in X$}
\end{align}
for some $T\in\mathcal{L}_+(Z,Y)$, a continuous linear operator with $T(z)\geq 0$ for all $z\in C$.
Without loss, $M\subseteq \dom f$.
Throughout this note, the operator $T\in\mathcal{L}_+(Z,Y)$  will take the special form $T=T_{(z^*,e)}$ for some $z^*\in C^+$ and $e\in\Int K$, defined as
\[
\forall z\in Z:\; T_{(z^*,e)}(z)=z^*(z)e
\]

\begin{definition}
An element $x_0$ is called an
\begin{enumerate}[(a)]
\item{\em $\eps$-v-$\wmin$-solution of \eqref{eq:P}}, iff $x_0\in M$ and there exists $y_0\in f(x_0)$ such that
\[
\of{\bigcup\limits_{x\in M}f(x)+\eps e}\cap\of{y_0+(-\Int K)}=\emptyset.
\]
A $0$-v-$\wmin$-solution is, for short, denoted as a v-$\wmin$-solution.
\item{\em $\eps$-l-$\wmin$-solution of \eqref{eq:P}}, iff $x_0\in M$ and for all $x\in M$ it holds
\[
f(x_0)\subseteq  \Int\of{f(x)+K}+\eps e\quad \Rightarrow\quad f(x_0)\subseteq \Int \of{f(x)+K}+\eps e .
\]
A $0$-l-$\wmin$-solution is, for short, denoted as a l-$\wmin$-solution.
\end{enumerate}
\end{definition}

The $v$ in the definition refers to the fact that $f(x_0)$ is $\eps$-weak minimal in the set $\cb{f(x)\st x\in M}$ with respect to the vector criteria, while the $l$ refers to the fact that $f(x_0)$ is $\eps$-weak minimal in the set $\cb{f(x)\st x\in M}$ with respect to the lattice criteria, compare \cite{FivePersonSoMF}.

\begin{lemma}\label{lem:solutionVSlsolution}
If $x_0\in M$ is a $\eps$-v-$\wmin$-solution to \eqref{eq:P}, then it is also a $\eps$-l-$\wmin$-solution to \eqref{eq:P}.
If additionally $f(x_0)=y_0+K$ is true for some $y_0\in f(x_0)$, then $x_0$ is a $\eps$-v-$\wmin$-solution to \eqref{eq:P}, if and only if it is a $\eps$-l-$\wmin$-solution to \eqref{eq:P}.
\end{lemma}
\proof
Let $x_0\in M$ be a $\eps$-v-$\wmin$-solution to \eqref{eq:P}, $y_0\in f(x_0)$ such that
\[
(y_0-\Int K)\cap \bigcup\limits_{x\in M}f(x)+\eps e=\emptyset
\]
and $f(x_0)\subseteq f(x)+\Int K+\eps e$, then $f(x)+\eps e\cap (y_0-\Int K)\neq\emptyset$, a contradiction.

On the other hand, let $f(x_0)=y_0+K$ be true for some $y_0\in f(x_0)$ and $y\in f(x)$ such that $x\in M$ and $y+\eps e\in (y_0-\Int K)$. Then $y_0\in f(x_0)\subseteq f(x)+\Int K+\eps e$ is true and as $x_0$ is a $\eps$-l-$\wmin$-solution to \eqref{eq:P}, this implies $y_0\in f(x)+\Int K+\eps e\subseteq f(x_0)+\Int K+2\eps e\subseteq y_0+\Int K$, a contradiction.
\pend

\begin{lemma}\label{lem:convexity}
Let $z^*\in C^+$ and $e\in\Int K$ be given.
If $\cl Q$ is convex, $z^*\in C^+$ and $e\in \Int K$, then
\[
\cl\of{\bigcup\limits_{x\in X}L(x,T_{(z^*,e)})+K}
\]
is convex, i.e.  if $(f\times g):X\to\P(Y\times Z)$ is closely convex like, then  $L(\cdot,T_{(z^*,e)}):X\to\P(Y)$ is closely convex like.
\end{lemma}
\proof
It holds
\[
\bigcup\limits_{x\in X}L(x,T_{(z^*,e)})+K=\cb{y+z^*(z)e\st (y,z)\in Q}
\]
and it is immediate that
\[
\cl\of{\bigcup\limits_{x\in X}L(x,T_{(z^*,e)})+K}\supseteq\cb{y+z^*(z)e\st (y,z)\in \cl Q}
\]
holds true. On the other hand, if $\bar y\notin\cl \of{\bigcup\limits_{x\in X}L(x,T_{(z^*,e)})+K}$, then it exists a convex $0$-neighbourhood $U\subseteq Y$ such that
\[
\forall (y,z)\in Q:\quad y+z^*(z)e\notin \bar y+U.
\]
Assume that $\bar y=y_0+z^*(z_0)e$ with $(y_0,z_0)\in\cl Q$, then there is a net in $Q$ converging to $(y_0,z_0)$ and  thus $y_i+z^*(z_i)e$ converges to $\bar y$, a contradiction, hence
\[
\cl\of{\bigcup\limits_{x\in X}L(x,T_{(z^*,e)})+K}=\cb{y+z^*(z)e\st (y,z)\in \cl Q},
\]
which immediately implies convexity.
\pend

In Langrange duality, the Slater condition plays an important role to guarantee strong duality. It already appears in the Lagrange multiplier rule. The following two equations are two possible generalizations to the set valued situation. We will prove in Lemma \ref{lem:Slater-CQ} that if $f\times g$ is closely convex like, then both are equivalent.

\begin{equation}\label{eq:Slater}
\exists x\in X:\quad g(x)\cap -\Int C\neq\emptyset 
\end{equation}

\begin{equation}\label{eq:CQ}
\forall z^*\in C^+\setminus\cb{0}\exists x\in\dom f\,\exists z\in g(x):\quad z^*(z)<0
\end{equation}

\begin{lemma}\label{lem:Slater-CQ}
Let $\cl Q$ be convex, then \eqref{eq:Slater} is satisfied, if and only if \eqref{eq:CQ} is satisfied.
\end{lemma}
\proof
\eqref{eq:Slater} immediately implies \eqref{eq:CQ}. On the other hand, let $\cl Q$ be convex and \eqref{eq:CQ} is satisfied, but $-\Int C\cap\bigcup\limits_{x\in X}g(x)=\emptyset$.
Thus for all $\bar y\in Y$ and $\bar c\in\Int C$, it holds $(\bar y,-\bar c)\notin\cl Q$. But as $\cl Q$ is convex, by a separation theorem (compare \cite[Theorem 1.1.6]{Zalinescu02} it holds
\[
\exists (y^*,z^*)\in K^+\times C^+\setminus\cb{(0,0)}:\quad y^*(\bar y)+z^*(-c)< \inf\limits_{(y,z)\in Q}y^*(y)+z^*(z).
\]
Let $\bar y$ be chosen such that $\bar y\in f(x)+\Int K=\Int (f(x)+K)$ is true for some $x\in X$, then
\[
\forall y^*\in K^+\setminus\cb{0}:\quad \inf\limits_{(y,z)\in Q}y^*(y)<y^*(\bar y),
\]
hence this choice of $\bar y$ implies $y^*=0$ and $z^*\neq 0$, implying
\begin{equation}\label{eq:1}
\forall c\in\Int C\,\exists z^* C^+\setminus\cb{0}:\quad z^*(-c)< \inf\limits_{(y,z)\in Q}z^*(z).
\end{equation}
For all $z^*\in C^+$ it holds
\[
\inf\limits_{(y,z)\in Q}z^*(z)=\inf\limits_{z\in \bigcup\limits_{x\in \dom f}g(x)+C}z^*(z),
\]
hence \eqref{eq:1} implies
\begin{equation*}
\forall c\in\Int C:\quad -c\notin \cl\co\of{\bigcup\limits_{x\in \dom f}g(x)+C},
\end{equation*}
hence $-\Int C\cap \cl\co\of{\bigcup\limits_{x\in \dom f}g(x)+C}=\emptyset$. But as both sets are convex, there exists a separating function $z^*\in C^+\setminus\cb{0}$, a contradiction to \eqref{eq:CQ}.
\pend

\section{Multiplier Rules in the Literature}

The result in \cite{corley87} is given in a form to maximize a concave function. The formulation below is equivalent to the original one, but given to fit the setting of this paper.

\begin{theorem}\label{MPR:Corley}\cite[Theorem 4.1]{corley87}
Let $f$ and $g$ be convex functions, \eqref{eq:Slater} is satisfied and $x_0\in M$ is a v-$\wmin$-solution to \eqref{eq:P}. Then there exists a $T\in\mathcal L_+(Z,Y)$ such that $x_0$ is a v-$\wmin$-solution to \eqref{eq:LP_T}.
\end{theorem}

In \cite{HernandezMarin2007lagrangian}, a collection of assumptions is given, such that a l-$\wmin$-solution to \eqref{eq:P} is also a l-$\wmin$-solution to \eqref{eq:LP_T}.

\begin{theorem}\label{MPR:HernMa}\cite[Theorem 4.1]{HernandezMarin2007lagrangian}
Let $Q$ be convex, \eqref{eq:Slater} is satisfied and $x_0\in M$ is a l-$\wmin$-solution to \eqref{eq:P}.
Define 
\[
B(x_0)=\cb{\bar y\in Y:\quad f(x_0)\subseteq \bar y+K};
\]
\[
H(x_0)=\cb{(y^*,z^*)\in K^+\setminus\cb{0}\times C^+\st 
\sup\limits_{y\in B(x_0)}(y^*,y)\leq \inf\limits_{(y,z)\in Q}(y^*,y)+(z^*,z)}
\]
\begin{enumerate}[(i)]
\item
$f(x_0)$ is $K$-bounded, that is
\item
$f(x_0)$ possesses $\Int K$-minimal elements, that is
\[
\exists y_0\in f(x_0):\quad f(x_0)\cap(y_0-\Int K)=\emptyset;
\]
\item
If $H\neq\emptyset$, then there is $(y^*_0,z^*_0)\in H$ and $(y_0,z_0)\in f(x_0)\times g(x_0)$ such that
\[
(y^*_0,y_0)+(z^*_0,z_0)= \inf\limits_{(y,z)\in Q}(y^*_0,y)+(z^*_0,z),
\]
\end{enumerate}
then there exists a $T\in\mathcal L_+(Z,Y)$ such that $x_0$ is a l-$\wmin$-solution to \eqref{eq:LP_T}. 
\end{theorem}

The formulation in \cite{HernandezMarin2007lagrangian} is given for $T+m$ with $m\in -K$, however the constant term does not influence the result and can be neglected. Moreover, the authors state that
and $T(z)+m\in -C$ for all $z\in g(x_0)\cap-C$. However this is immediate, as $T\in\mathcal{L}_+(Z,Y)$.

Most recently, a multiplier rule for approximate solutions is given in \cite{DhingraLalitha2016}.
\begin{theorem}\label{MPR:Dhing}\cite[Theorem 3.2]{DhingraLalitha2016}
Let $\cl Q$ be convex, \eqref{eq:CQ} is satisfied and $x_0\in M$ is a $\eps$-l-$\wmin$-solution to \eqref{eq:P}.
If additionally 
\[
\exists y_0\in f(x_0):\quad f(x_0)\subseteq y_0+K,
\]
then there exists a $T\in\mathcal L_+(Z,Y)$ such that $x_0$ is a $\eps$-l-$\wmin$-solution to \eqref{eq:LP_T}  and
\[
\forall z\in g(x_0)\cap-C:\quad T(z)=0.
\]
\end{theorem}

All three theorems are proven by setting $T=T_{(z^*,e)}$ for a specified $z^*\in C^+$, $e\in \Int K$. It is easily seen that the convexity assumption in \cite{DhingraLalitha2016} is the most general one. However, all three theorems are proven by using a separation theorem which only applies the convexity of $\cl Q$.

In the next section we will prove that, but for the convexity assumption, the assumptions of Theorem \ref{MPR:Corley} are the most general ones of the three and those of Theorem \ref{MPR:Dhing} the most specific for the case $\eps=0$. 
Theorem \ref{thm:Solution} will provide an argument why the assumptions of Theorem \ref{MPR:Corley} can only be weakened in terms of the convexity assumption and Theorem \ref{thm:eps-Solution} generalizes Theorem \ref{thm:Solution} to $\eps$-v-$\wmin$-solutions. This result is also generalization of the result given in \cite{DhingraLalitha2016}.

\section{Connecting the Dots}

\begin{lemma}\label{lem:eps-P-Solution}
Let $\cl Q$ be convex, then  $x_0\in X$ is a $\eps$-v-$\wmin$-solution to \eqref{eq:P}, if and only if
\[
\exists y_0\in f(x_0)\;
\exists (y^*,z^*)\in K^+\times C^+\setminus\cb{(0,0)}:\quad \; 
y^*(y_0-\eps e)\leq\inf\limits_{(y,z)\in Q}y^*(y)+z^*(z).
\]
In this case, $-\eps y^*(e)\leq z^*(z)$ is true for all $z\in g(x_0)$.

If additionally \eqref{eq:Slater} is satisfied, then $y^*\neq 0$ is true and without loss $y^*(e)=1$.

\end{lemma}
\proof
Assume $x_0\in X$ is a $\eps$-v-$\wmin$-solution to \eqref{eq:P}, i.e.
\[
\exists y_0\in f(x_0):\quad \of{y_0-\eps e-\Int K}\cap \bigcup\limits_{x\in M}f(x)=\emptyset
\] 
and $(y,z)\in Q\cap\of{y_0-\eps e-\Int K\times-\Int C}$. Then by definition there is $x\in M$ such that $y\in f(x)$ and $z\in g(x)$, hence
\[
y\in \of{y_0-\eps e-\Int K}\cap \bigcup\limits_{x\in M}f(x),
\]
a contradiction. Thus, if $x_0$ is a $\eps$-v-$\wmin$-solution to \eqref{eq:P}, then
\[
\exists y_0\in f(x_0):\quad \of{y_0-\eps e-\Int K\times-\Int C}\cap \cl Q=\emptyset.
\]
Both sets are convex, hence 
\begin{align*}
\exists (y^*,z^*)&\in K^+\times C^+\setminus\cb{(0,0)}:\\
&\sup\limits_{(y,z)\in \Int \of{K\times C}}y^*(y_0-y-\eps e)+z^*(-z)\leq \inf\limits_{(y,z)\in Q}y^*(y)+z^*(z),
\end{align*}
which is equivalent to 
\[
\exists (y^*,z^*)\in K^+\times C^+\setminus\cb{(0,0)}:\quad 
y^*(y_0-\eps e)\leq\inf\limits_{(y,z)\in Q}y^*(y)+z^*(z).
\]
But as
\[
\inf\limits_{(y,z)\in Q}y^*(y)+z^*(z)\leq \inf\limits_{\substack{ y\in f(x), z\in g(x),\\x\in M}}y^*(y)+z^*(z)
\leq \inf\limits_{\substack{ y\in f(x),\\x\in M}}y^*(y)
\]
is true for all $(y^*,z^*)\in K^+\times C^+\setminus\cb{(0,0)}$, this in turn implies
\[
(y_0-\eps e-\Int K)\cap \co\bigcup\limits_{x\in M}f(x)=\emptyset,
\] 
hence $x_0$ is a $\eps$-v-$\wmin$-solution to \eqref{eq:P}.
Moreover,
\[
y^*(y_0-\eps e)\leq\inf\limits_{(y,z)\in Q}y^*(y)+z^*(z)\leq \inf\limits_{(y,z)\in f(x_0)\times g(x_0)}y^*(y)+z^*(z)
\]
implies $-\eps y^*(e)\leq \inf\limits_{z\in g(x_0)}z^*(z)$.

Finally, let $g(x)\cap-\Int C\neq\emptyset$ be assumed, $y_0\in f(x_0)$, $y^*=0$ and $z^*\in C^+$ be such that
\[
y^*(y_0-\eps y^*(e))\leq\inf\limits_{(y,z)\in Q}y^*(y)+z^*(z),
\]
that is
\[
0\leq \inf\limits_{(y,z)\in Q}z^*(z),
\]
a contradiction to \eqref{eq:Slater}.

\pend

An alternative formulation of Lemma \ref{lem:eps-P-Solution} is as follows

\begin{lemma}\label{lem:P-Solution2}
Let $\cl Q$ be convex, then  $x_0\in X$ is a $\eps$-v-$\wmin$-solution to \eqref{eq:P}, if and only if
\begin{align*}
\exists (y_0,z_0)\in f(x_0)\times g(x_0)\,
&\exists (y^*,z^*)\in K^+\times C^+\setminus\cb{(0,0)}:\\ 
&y^*(y_0-\eps e)+ z^*(z_0)\leq\inf\limits_{(y,z)\in Q}y^*(y)+z^*(z).
\end{align*}
In this case, $-\eps y^*(e)\leq z^*(z)$ is true for all $z\in g(x_0)$.

If additionally \eqref{eq:Slater} is satisfied, then $y^*\neq 0$ is true.
\end{lemma}

\begin{corollary}
Let $\cl Q$ be convex, then  $x_0\in X$ is a v-$\wmin$-solution to \eqref{eq:P}, if and only if
\begin{align*}
\exists (y_0,z_0)\in f(x_0)\times g(x_0)\,
&\exists (y^*,z^*)\in K^+\times C^+\setminus\cb{(0,0)}:\\ 
&y^*(y_0)+ z^*(z_0)=\inf\limits_{(y,z)\in Q}y^*(y)+z^*(z).
\end{align*}
In this case, $z^*(z)=0$ is true for all $z\in g(x_0)$.

If additionally \eqref{eq:Slater} is satisfied, then $y^*\neq 0$ is true.
\end{corollary}

\begin{lemma}\label{lem:eps-LP-Solution}
Let $\cl Q$ be convex, $z^*\in C^+$ and $e\in \Int K$, then  $x_0\in X$ is a $\eps$-v-$\wmin$-solution to \eqref{eq:LP_T} with $T=T_{(z^*,e)}$, if and only if
\begin{align*}
\exists (y_0,z_0)\in f(x_0)\times g(x_0)\,
&\exists y^*\in K^+\setminus\cb{0},\, y^*(e)=1,\\
&y^*(y_0)+z^*(z_0)-\eps\leq\inf\limits_{(y,z)\in Q}y^*(y)+z^*(z).
\end{align*}
\end{lemma}
\proof
By definition $x_0\in X$ is a $\eps$-v-$\wmin$-solution to \eqref{eq:LP_T} with $T=T_{(z^*,e)}$, if and only if
\[
\exists (y_0,z_0)\in f(x_0)\times g(x_0):\quad 
(y_0+z^*(z_0)e)-\eps e\notin\Int \cb{y+z^*(z)e\st (y,z)\in Q}.
\]
As by Lemma \ref{lem:convexity} $\cl\cb{y+z^*(z)e\st (y,z)\in Q}$ is convex, this is equivalent to
\begin{align*}
\exists (y_0,z_0)\in f(x_0)\times g(x_0)\,
&\exists y^*\in K^+\setminus\cb{0},\, y^*(e)=1,\\
&y^*(y_0)+z^*(z_0)-\eps\leq\inf\limits_{(y,z)\in Q}y^*(y)+z^*(z).
\end{align*}

\pend

\begin{corollary}
Let $\cl Q$ be convex, $z^*\in C^+$ and $e\in \Int K$, then  $x_0\in X$ is a v-$\wmin$-solution to \eqref{eq:LP_T} with $T=T_{(z^*,e)}$, if and only if
\begin{align*}
\exists (y_0,z_0)\in f(x_0)\times g(x_0)\,
&\exists y^*\in K^+\setminus\cb{0},\, y^*(e)=1,\\
&y^*(y_0)+z^*(z_0)=\inf\limits_{(y,z)\in Q}y^*(y)+z^*(z).
\end{align*}
\end{corollary}

\begin{theorem}\label{thm:eps-Solution}
Let $\cl Q$ be convex and \eqref{eq:Slater} be satisfied. Then  $x_0\in X$ is a $\eps$-v-$\wmin$-solution to \eqref{eq:P}, if and only if it exists $z^*\in C^+$ and $e\in \Int K$ such that $x_0$ is a $\eps$-v-$\wmin$-solution to \eqref{eq:LP_T} with $T=T_{(z^*,e)}$. In this case, $T(g(x_0))\subseteq\cb{(t-\eps)e\st t\geq 0}$.
\end{theorem}
\proof
Proven in Lemma \ref{lem:eps-P-Solution} and Lemma \ref{lem:eps-LP-Solution}.
\pend

\begin{remark}
Under the assumption of \cite[Theorem 3.2]{DhingraLalitha2016}, $\cl Q$ is convex, \eqref{eq:Slater} is satisfied and  $x_0\in X$ is a $\eps$-v-$\wmin$-solution to \eqref{eq:P}, hence it is a special case of Theorem \ref{thm:eps-Solution}, while the conclusions of Theorem \ref{thm:eps-Solution} is more specific than those in \cite{DhingraLalitha2016}, hence our result strengthens the previous one.
\end{remark}

\begin{theorem}\label{thm:Solution}
Let $\cl Q$ be convex and \eqref{eq:Slater} is satisfied, then  $x_0\in X$ is a v-$\wmin$-solution to \eqref{eq:P}, if and only if it exists $z^*\in C^+$ and $e\in \Int K$ such that $x_0$ is a v-$\wmin$-solution to \eqref{eq:LP_T} with $T=T_{(z^*,e)}$. In this case, $T(z)=0$ for all $z\in g(x_0)\cap -C$.
\end{theorem}

\begin{remark}
The assumptions in Theorem \ref{thm:Solution} differ from those in \cite[Theorem 4.1]{corley87} only in a relaxed convexity assumption on $f\times g$, while those in \cite[Theorem 3.2]{DhingraLalitha2016} (for $\eps=0$) are more specific than those of Theorem \ref{thm:Solution}, compare Lemma \ref{lem:solutionVSlsolution}. 

In \cite[Theorem 4.1]{corley87} it is proven that the assumptions imply the existence of a $T\in \mathcal{L}_+(Z,Y)$ such that $x_0$ is a v-$\wmin$-solution to \eqref{eq:LP_T} while in \cite[Theorem 3.2]{DhingraLalitha2016} the existence of  a $T\in \mathcal{L}_+(Z,Y)$ such that $x_0$ is a l-$\wmin$-solution to \eqref{eq:LP_T} and $T(z)=0$ for all $g(x_0)\cap-C$ is proven. As any v-$\wmin$-solution is especially a l-$\wmin$-solution, the result in \cite[Theorem 4.1]{corley87} is slightly more specific in this instance while Theorem \ref{thm:Solution} is more specific than both.
\end{remark}

\begin{lemma}\label{lem:l-P-Solution}
Let $\cl Q$ be convex, $x_0\in M$
 and $B=\cb{\tilde y\in Y\st f(x_0)\subseteq \tilde y+K}$. If  $x_0\in X$ is a l-$\wmin$-solution to \eqref{eq:P}, then
\[
\exists (y^*,z^*)\in K^+\times C^+\setminus\cb{(0,0)}:\quad 
\sup\limits_{\tilde y\in B}y^*(\tilde y)\leq\inf\limits_{(y,z)\in Q}y^*(y)+z^*(z).
\]
If additionally \eqref{eq:Slater} is assumed, $y^*\in K^+\setminus\cb{0}$ is true.
\end{lemma}
\proof
Let $x_0$ be a l-$\wmin$-solution to \eqref{eq:P} and
assume to the contrary
\[
\exists (\bar y,\bar c)\in Y\times C\,\exists \tilde y\in B:\quad (\bar y,-\bar c)\in (\tilde y-\Int K\times-\Int C)\cap Q,
\]
hence there exists $\bar x\in M$ such that 
\[
f(x_0)\subseteq\tilde y+K\subseteq f(\bar x)+\Int K.
\]
 As $x_0$ is a l-$\wmin$-solution to \eqref{eq:P}, this implies 
\[
\tilde y+K\subseteq f(\bar x)+K\subseteq f(x_0)+\Int K\subseteq \tilde y+\Int K,
\]
a contradiction.
Thus 
\[
(\Int B\times-\Int C)\cap \cl Q=\emptyset
\]
and as both sets are convex, they can be separated by a continuous linear function $(y^*,z^*)\in K^+\times C^+\setminus\cb{(0,0)}$, implying
\[
\exists (y^*,z^*)\in K^+\times C^+\setminus\cb{(0,0)}:\quad \; 
\sup\limits_{\tilde y\in B}y^*(\tilde y)\leq\inf\limits_{(y,z)\in Q}y^*(y)+z^*(z).
\]
Under the assumption of \eqref{eq:Slater} and $y^*=0$ this implies
\[
0\leq\inf\limits_{(y,z)\in Q}z^*(z),
\]
a contradiction, hence $y^*\in K^+\setminus\cb{0}$.
\pend

\begin{remark}
By Lemma \ref{lem:l-P-Solution}, under the assumptions of \cite[Theorem 4.1]{HernandezMarin2007lagrangian} it holds
\begin{align*}
\exists (y_0,z_0)\in f(x_0)\times g(x_0)\,
&\exists y^*\in K^+\setminus\cb{0},\, y^*(e)=1,\\
&y^*(y_0)+ z^*(z_0)=\inf\limits_{(y,z)\in Q}y^*(y)+z^*(z).
\end{align*}
But by Lemma \ref{lem:P-Solution2} this is true, if and only if $x_0$ is a v-$\wmin$-solution to \eqref{eq:P}, thus \cite[Theorem 4.1]{HernandezMarin2007lagrangian} is a special case of Theorem \ref{thm:Solution}.

On the other hand, under the assumptions of \cite[Theorem 3.2]{DhingraLalitha2016} (for $\eps=0$), $f(x_0)$ is bounded by $y_0\in f(x_0)$, $x_0$ is a v-$\wmin$-solution to \eqref{eq:P}, $B=y_0-K$ and by Lemma \ref{lem:P-Solution2}
\begin{align*}
\exists (y_0,z_0)\in f(x_0)\times g(x_0)\,
&\exists y^*\in K^+\setminus\cb{0},\, y^*(e)=1,\\
&y^*(y_0)+z^*(z_0)=\inf\limits_{(y,z)\in Q}y^*(y)+z^*(z).
\end{align*}
Hence the assumptions are more specific, but for the convexity assumption, than those of \cite[Theorem 4.1]{HernandezMarin2007lagrangian}, while the implication in \cite[Theorem 4.1]{HernandezMarin2007lagrangian} is more specific on the fact, that $x_0$ is a v-$\wmin$-solution to \eqref{eq:LP_T}.

\end{remark}

\section{Conclusion}
In Theorem \ref{thm:eps-Solution} and \ref{thm:Solution}, a Lagrange multiplier rule for (approximate) weak solutions (with respect to the vector criteria) to a constrained set valued optimization problem is given. Both rules provide an equivalent description of a (approximate) weak solution. We have proven that, when relaxing the convexity assumption to $\cl Q$ convex in the older two results, the multiplier rule in \cite{DhingraLalitha2016} is, for $\eps=0$, a special case of the one in \cite{HernandezMarin2007lagrangian}, which in turn is a special case of the one presented in \cite{corley87}. Ultimately, the assumptions in Theorem \ref{thm:Solution} are identical to those in \cite{corley87} but for the convexity assumption, so our result is a generalization of all three multiplier rules under consideration. Notably, Theorem \ref{thm:eps-Solution} is a generalization of the multiplier rule for approximate solutions presented in \cite{DhingraLalitha2016}.
The implications in Theorem  \ref{thm:Solution} are given in a stronger way than in the quoted multiplier rules. In the proofs for all three multiplier rules, the authors exploit the fact that the implication is of the form we stated in  Theorem \ref{thm:Solution}. Thus our conclusion is stronger, but implicitly present in the quoted literature.

\end{document}